\documentclass[11pt]{article}
\usepackage{amsmath,amssymb,amscd,array}

\let\ssection=\section
\renewcommand{\section}{\setcounter{equation}{0}\ssection}

\def\d{\delta}

\def\om{\omega}
\def\r{\rho}
\def\a{\alpha}
\def\b{\beta}
\def\s{\sigma}
\def\vfi{\varphi}

\def\l{\lambda}
\def\m{\mu}

\def\implies{\Rightarrow}

\newcommand{\Diff}{\mathrm{Diff}}

\newcommand{\cF}{{\mathcal{F}}}

\newcommand{\Hom}{\mathrm{Hom}}
\newcommand{\ccD}{\mathcal{D}}

\newcommand{\Vect}{\mathrm{Vect}}

\newcommand{\cqfd}{\hspace*{\fill}\rule{3mm}{3mm}}
\newcommand{\cqf}{\hspace*{\fill}\rule{2mm}{2mm}}

\begin{document}

\frenchspacing

\def\d{\delta}
\def\g{\gamma}
\def\om{\omega}
\def\r{\rho}
\def\a{\alpha}
\def\b{\beta}
\def\s{\sigma}
\def\vfi{\varphi}
\def\l{\lambda}
\def\m{\mu}
\def\implies{\Rightarrow}

\oddsidemargin .1truein
\newtheorem{thm}{Theorem}[section]
\newtheorem{lem}[thm]{Lemma}
\newtheorem{cor}[thm]{Corollary}
\newtheorem{pro}[thm]{Proposition}
\newtheorem{ex}[thm]{Example}
\newtheorem{rmk}[thm]{Remark}
\newtheorem{defi}[thm]{Definition}
\title{The space of $m$-ary differential operators as a module
over the Lie algebra of vector fields }
\author{Sofiane BOUARROUDJ\footnote{Supported by the
Grant 03130211/05, U.A.E. University}\\
{\footnotesize Department of Mathematical Sciences, U.A.E.
University} \\
{\footnotesize P.O. Box 17551, Al-Ain, U.A.E.}\\
{\footnotesize  e-mail:bouarroudj.sofiane@uaeu.ac.ae} }
\date{}
\maketitle
\begin{abstract}
The space $\ccD_{\underline{\lambda};\mu},$ where
$\underline{\lambda}=(\l_1,\ldots,\l_m),$ of $m$-ary differential
operators acting on weighted densities is a $(m+1)$-parameter
family of modules over the Lie algebra of vector fields. For
almost all the parameters, we construct a canonical isomorphism
between the space $\ccD_{\underline{\lambda};\mu}$ and the
corresponding space of symbols as $\mathfrak{sl}(2)$-modules. This
yields to the notion of the $\mathfrak{sl}(2)$-equivariant symbol
calculus for $m$-ary differential operators. We show, however,
that these two modules cannot be isomorphic as
$\mathfrak{sl}(2)$-modules for some particular values of the
parameters. Furthermore, we use the symbol map to show that all
modules $\ccD_{\underline{\lambda};\mu}^2$ (i.e., the space of
second-order operators) are isomorphic to each other, except for
few modules called singular.
\end{abstract}
\section{Introduction}
In this paper let $M$ be either $\mathbb{R}$ or $S^1;$ let
$\cF_\l$ be the space of weighted densities on $M$ of weight $\l,$
i.e., the space of sections of the line bundle $ (T^*M)^{\otimes
\lambda}$, where $\lambda\in \mathbb{R}$. This space has the
following structure as a $\mathrm{Vect}(M)$-module: for any
$a\,(dx)^{\lambda}\in {\cal F}_{\lambda}$ and $X\in
\mathrm{Vect}(M),$ we put
\begin{equation}
\label{dens} L_{X}^{\lambda}(a\,(dx)^{\lambda})=(X ({a}) +\lambda
\,a \,\text{div}\,X)\,(dx)^{\lambda}.
\end{equation} Denote
by $\ccD_{\underline{\l};\mu},$ where
$\underline{\lambda}=(\l_1,\ldots,\l_m),$ the space of $m$-ary
linear differential operators:
\[
\cF_{\l_1}\otimes \cdots \otimes \cF_{\l_m} \rightarrow \cF_\mu.
\]
The action by the Lie derivative on $\ccD_{\underline{\l};\mu}$
defines a module structure over the Lie algebra of vector fields,
$\Vect(M).$ The space of unary differential operators viewed as a
$\Vect(M)$-module is a classical object (see, e.g., \cite{wil}).

In this paper, we study the geometry of the modules
$\ccD_{\underline{\l};\mu}$ understood in the sense of Klein: as Lie
group actions (or Lie algebra actions) on a manifold. We will be
dealing with the Lie algebra of smooth vector fields, $\Vect(M),$
and $\mathfrak{sl}(2)$ embedded into $\Vect(M)$ via infinitesimal
projective transformations:
\begin{equation}
\label{Dub} \mathfrak{sl}(2) \simeq\mathrm{Span}\{\partial_x,
x\,\partial_x, x^2\partial_x\}.
\end{equation}
The quotient module
$\ccD_{\underline{\l};\mu}^k/\ccD_{\underline{\l};\mu}^{k-1}$ can
be decomposed into $\binom{k+m-1}{m-1}$ components that transform
under coordinates change as $(\delta-k)$-densities, where
$\delta=\mu-\sum_{j=1}^m\l_j.$ Therefore, the multiplication of
these components by any non-singular matrix, say $\alpha,$ gives
rise to an isomorphism
\[
\sigma^{\alpha}:
\ccD_{\underline{\l};\mu}^k/\ccD_{\underline{\l};\mu}^{k-1}
\xrightarrow{\simeq}
\cF_{\delta-k} \oplus\cdots \oplus\cF_{\delta-k}.
\]
The map $\sigma^{\alpha}$ is what we call the {\it principal
symbol}. By the very definition, the principal symbol is
$\Vect(M)$-equivariant but not unique if $m>1.$ Let us consider
the graded space
\[
{\cal S}_\delta= \bigoplus_{k\geq 0} {\cal
D}^k_{\underline{\lambda};\mu}/{\cal
D}^{k-1}_{\underline{\lambda};\mu}
\]
associated with the natural filtration of ${\cal
D}_{\underline{\lambda};\mu}.$ A {\it symbol} map is a linear
bijection
\[
\sigma^{\alpha}_{\underline{\l},\mu}:
\ccD_{\underline{\l};\mu}\longrightarrow  {\cal S}_\delta
\]
such that the highest-order term of
$\sigma^{\alpha}_{\underline{\l},\mu}(A),$ where $A\in
\ccD_{\underline{\l};\mu},$ coincides with the principal symbol
$\sigma^\alpha(A).$ In the unary case, equivariant symbol calculus
has been first introduced in \cite{cmz}, then studied in \cite{gar}.
A generalization to multi-dimensional manifolds has been first
studied in \cite{dlo, lo1}, then studied in
\cite{bor,bl,bhm,hans,loub,mr}. In the binary case, the existence of
the symbol map has been investigated in \cite{bn} and an explicit
formula has been given for the space of second-order operators. We
show, for almost all $\underline{\l}$ and $\mu,$ that there exists a
$\mathfrak{sl}(2)$-equivariant symbol map for every $m.$ As in the
unary case, we show that the symbol maps are not
$\Vect(M)$-equivariant except for $k=1$ or $k=2$ but for particular
values of $\underline{\l}$ and $\mu.$ Although the symbol map is not
unique for $m>1$ -- unlike the unary case -- the uniqueness can be
understood as follows: once the the symbol map $\sigma^\alpha$ is
fixed, the corresponding $\mathfrak{sl}(2)$-equivariant map
$\sigma^{\alpha}_{\underline{\l},\mu}$ is unique.

Furthermore, we investigate for which parameters
$(\underline{\l},\mu)$ and $(\underline{\rho},\eta),$ we have
$\ccD_{\underline{\l};\mu}^k~\simeq~
\ccD_{\underline{\rho};\eta}^k$ as $\Vect(M)$-modules. We only
deal with $k=1,2.$ In the unary case, the classification problem
of such modules has first been raised and studied in \cite{do};
and comprehensive results were obtained in the papers \cite{do,
gar, lmt}. As for the unary case, a necessary condition for the
isomorphism is imposed by
\[
\delta:=\mu-\sum_{j=1}^m\l_j=\eta-\sum_{j=1}^m\rho_j.
\]
We prove that all modules are isomorphic to each other, provided
they have the same $\delta,$ except:
\begin{itemize}
    \item the module ${\cal D}^1_{0,\ldots,0;\delta};$
    \item the modules ${\cal D}^2_{0,\ldots,0;\delta}\simeq
    {\cal D}^2_{1-\delta,\ldots,0;1}\simeq \ldots \simeq
    {\cal D}^2_{0,\ldots,1-\delta;1}.$
\end{itemize}
In the generic case, we use the symbol map
$\sigma^{\alpha}_{\underline{\l},\mu}$ to build-up the isomorphism
between $\ccD_{\underline{\l};\mu}^k\quad \mbox{ and }\quad
\ccD_{\underline{\rho};\eta}^k.$ However, for some values of
$\delta,$ called {\it resonant}, we proceed by a direct computation.
The proof in all cases is based on the locality of any
$\Vect(M)$-isomorphism.
\section{The space of $m$-ary differential operators as
a $\Vect(M)$-module}
We fix a natural number $m.$ In order to avoid clutter, we have
found that it is convenient to use the following notations:
\begin{itemize}
    \item Denote by $\underline{i}$ either the $m$-tuple
$(i_1,\ldots,i_m)$ or the indices $i_1,\ldots,i_m,$ as, for
instance, $a_{\underline{i}}=a_{i_1,\ldots,i_m}.$ The difference
should be discernable from the context.
     \item Denote by $|\underline{i}|$ the sum $\sum_{j=1}^mi_j.$
     \item Denote by $[\tau]_p$ the
     square matrix $[\tau^{\underline{i}}_{\underline{j}}]$ of
     size $\binom{p+m-1}{m-1}\times \binom{p+m-1}{m-1}$ for
     $|\underline{i}|=|\underline{j}|=p.$
    \item Denote $\mathbf{1}_i:=(0,\ldots,0,1,0\ldots,0),$ where 1 is in
    the $i$-th position.
    \item Denote by $ {\cal S}_{\lambda}^{(i)}=\bigoplus \cF_\l,$ where
$\cF_\lambda$ is counted $ {\binom{i+m-1}{m-1} \mbox{ times. }}$
\end{itemize}

Consider $m$-ary differential operators that act on weighted
densities:
\begin{equation}
\label{Op} A:{\cal F}_{\lambda_1}\otimes \cdots \otimes {\cal
F}_{\lambda_m}\to {\cal F}_{\mu}\quad \underline{\varphi} \mapsto
\sum_{r\geq 0}\sum_{|\underline{i}|=r}a_{\underline{i}}\,
\partial_{\underline{i}}\;\underline{\varphi},
\end{equation}
where $a_{\underline{i}}$ are smooth functions on $M.$ We denote by
${\cal D}_{\underline{\lambda};\mu}^k,$ where
$\underline{\lambda}=(\l_1,\ldots,\l_m),$ the space of $k$th-order
$m$-ary differential operators (\ref{Op}) endowed with the following
$\Vect(M)$-module structure.

For all $X\in \Vect(M),$ we define ($L_{X}^{\l_j}$ is the action
(\ref{dens})):
\begin{equation}
\label{act} L_{X}^{\underline{\lambda};\mu}(A)=L_{X}^\mu \circ A
-\sum_{j=1}^{m}A(...,L_{X}^{\l_{j}}(-),...).
\end{equation}
Let us denote by $a_{\underline{i}}^X$ the coefficients of the
operator $L_{X}^{\underline{\lambda};\mu}(A).$ We have
\begin{pro} \label{prop21}
The coefficients $a_{\underline{i}}^X$ can be expressed in terms of
the coefficients $a_{\underline{i}}$ as follows (where
$\delta=\mu-|\underline{\lambda}|$):
\begin{equation}
\label{actm}
\small
\begin{array}{ccl}
a^X_{\underline{s}}&=&L_{X}^{\delta-|\underline{s}|}
a_{\underline{s}}\\
&&-\displaystyle\sum_{i\geq s_1+1}^k {\small
\binom{i}{i+1-s_1}}\,X^{(i+1-s_1)}\,
a_{\underline{s}_{|_{s_1=i}}}-\cdots-\sum_{i\geq s_m+1}^k
\binom{i}{i+1-s_m}\,X^{(i+1-s_m)}
\,a_{\underline{s}_{|_{s_m=i}}}\\
&& -\displaystyle  \l_1 \sum_{i\geq s_1+1}^k
\binom{i}{i-s_1}\,X^{(i+1-s_1)} \,a_{\underline{s}_{|_{s_1=i}}}-
\cdots -\l_m\sum_{i\geq s_m+1}^k \binom{i}{i-s_m}\,X^{(i+1-s_m)}
\,a_{\underline{s}_{|_{s_m=i}}},\\
\displaystyle a^X_{\underline{0}}&=&\displaystyle
L_X^{\delta}a_{\underline{0}}-\l_1 \sum_{i=1}^{k}
X^{i+1}\,a_{i\mathbf{1}_1} -\cdots -\l_m\sum_{i=1}^{k} X^{i+1}
\,a_{i\mathbf{1}_m}.\\
\end{array}
\end{equation}
\end{pro}
{\bf Proof.} These formulas come out easily from the definition
(\ref{act}). \cqf

\section{Locality of the diffeomorphism, the invariant $\delta$}
A map $T$ is called {\it local} if $\mathrm{Supp} (T(A))\subset
\mathrm{Supp}(A)$ for all $A\in \ccD_{\underline{\l};\mu}.$ The
well-known Peetre's theorem (cf. \cite{Pee}) asserts that such a
map is a differential operator. Not only is $T(A)$ a differential
operator that acts on weighted densities but its coefficients are
given by a differential operator as well. The following
proposition is adapted from the unary case (see \cite{lo1}).
\begin{pro}
\label{main6} For $k\leq 2,$ every $\Vect(M)$-equivariant
isomorphism $T:\ccD_{\underline{\l};\mu}^k\rightarrow
\ccD_{\underline{\rho};\eta}^k$ is local.
\end{pro}
{\bf Proof.} Assume that $A\in\ccD_{\underline{\l};\mu}^k$ vanish
on an open subset $U\subset M.$ We will show that $T(A)$ vanishes
on $U$ as well. We have two cases:\\

1) The case where $\eta-|\underline{\rho}|\not = 0,1,2.$ Suppose
the contrary, namely $T(A)_{|_u}\not =0$ for some $u\in U.$ The
principle symbol $\sigma(T(A))$ of $T(A)$ can be expressed as
$\oplus_{i} \sigma_i\,(dx)^{\eta-|\underline{\rho}|-p},$ for a
certain integer $p.$ We can always choose $u$ such that
${\sigma_{i_0}}_{_{|_u}}\not =0,$ for a certain $i_0.$ To get the
contradiction we will look for a vector field $X$ that satisfies
$L_X^{\underline{\rho};\eta}(T(A))_{|_u}\not =0.$ Hence the
contradiction since
\[
L_X^{\underline{\lambda}; \eta} (A_{|_u})=0 \quad \mbox{and} \quad
T\left ( L_X^{\underline{\lambda}; \eta} (A_{|_u})\right
)=L_X^{\underline{\rho}; \eta} (T(A)_{|_u}).
\]
To choose $X$ we consider the expression
\begin{equation}
\label{we} L_X^{\eta- |\underline{\rho}|-p}( \sigma_{i_0}
(dx)^{\eta- |\underline{\rho}|-p})=\left (X(\sigma_{i_0}) +(\eta
-|\underline{\rho}|-p)\;\mathrm {div}\;X\sigma_{i_0}\right )
(dx)^{\eta- |\underline{\rho}|-p}.
\end{equation}
Since $\eta-|\underline{\rho}|\not = 0,1,2$ and $p\leq 2,$ we can
always choose $X$ such that Eq. (\ref{we}) is not zero. This
implies that $L_X^{\underline{\rho};\eta}(T(A))_{|_u}\not =0$ and
hence the contradiction.\\

2) The case where $\eta-|\underline{\rho}| = 0,1,2.$ Let
$\Vect_U(M)$ be the Lie algebra of smooth vector fields with
support in $U.$ Since $A$ vanishes on $U,$ it follows that
\[
L_X^{\underline{\lambda}; \mu}(A)=0\quad \mbox{ for every } X\in
\Vect_U(M).
\]
Therefore, $L_X^{\underline{\rho}; \eta}(T(A))=0$ for every $X\in
\Vect_U(M).$ This means that the operator $T(A)$ is
$\Vect_U(M)$-invariant, as $T(A)$ does not vanish on $U.$
Following \cite{groz,kir}, such an operator can be expressed as
follows:

(i) $A$ is the multiplication operator for
$\eta-|\underline{\rho}|=0$,

(ii) $A=\sum_{i,j=1}^m a_{i,j} \{\cdot,\cdot\},$ where $a_{i,j}$
are scalars and $\{\cdot , \cdot \}$ is the Poisson bracket, for
$\eta-|\underline{\rho}|=1$,

(iii) $A$ is a linear combination of operators given by
compositions of the de Rham operator and the Poisson bracket for
$\eta-|\underline{\rho}|=2$ and for special values of $\rho$ and
$\underline{\eta}.$\\
Now any operator as in (i)--(iii) is not only
$\Vect_U(M)$-invariant but also $\Vect(M)$-invariant. The
isomorphism $T$ implies that the operator $A$ is
$\Vect(M)$-invariant and is given as in (i)--(iii). This is a
contradiction, since if $A$ vanishes on $U$ it must vanish
everywhere. \cqf
\begin{pro}
\label{mai} Every linear $\mathfrak{sl}(2)$-equivariant
isomorphism $\ccD_{\underline{\l};\mu} \rightarrow {\cal
S}_{\delta} $ is local.
\end{pro}
{\bf Proof.} The statement of this proposition, in fact, holds
also true for the affine Lie sub-algebra ${\mathfrak a}=\mbox{Span
}\{\partial_x, x\partial_x \}$ of $\mathfrak{sl}(2).$ As
${\mathfrak a}$-modules, the space $\ccD_{\underline{\l};\mu}$ and
the space ${\cal S}_\delta$ are isomorphic, thanks to Proposition
\ref{prop21}. Therefore, it is sufficient to prove that any
${\mathfrak a}$-equivariant linear map
$\cF_{\delta-|\underline{i}|}\rightarrow
\cF_{\delta-|\underline{j}|}$ is local. This has been proved in
\cite{lo1} (Theorem  5.1) using Petree's theorem for $\delta=0$
but the proof works well for every $\delta.$  \cqf\\

For every module $\ccD_{\underline{\l};\mu}^k,$ we define its {\it
shift} $\delta$ to be
\[
\delta:=\mu-|\underline{\l}|.
\]
\begin{pro}
A necessary condition for the two modules
$\ccD_{\underline{\l};\mu}^k$ and $\ccD_{\underline{\rho};\eta}^k$
to be isomorphic is to have the same shift.
\end{pro}
{\bf Proof.} Let $T:\ccD_{\underline{\l};\mu}^k \rightarrow
\ccD_{\underline{\rho};\eta}^k$ be a $\Vect(M)$-isomorphism. We
shall study the equivariance property with respect to the vector
field $X=x\partial_x$ upon taking $\underline{\phi}$ is constant.
The $\Vect(M)$-equivariance reads as follows:
\[
\left (T(L_{X}^{\underline{\l};\mu} A)\right
)(\underline{\phi})=\left (L_{X}^{\underline{\rho};\eta} T(A)\right
)(\underline{\phi}).
\]
Consider $A=a_{\underline{0}} \partial_{\underline{0}},$ where
$a_{\underline{0}}$ is a smooth function on $M,$ the operator of
multiplication. By using Proposition \ref{prop21}, we have $L_X
^{\underline{\l};\mu}A=
(L_X^{\mu-|\underline{\l}|}\,a_{\underline{0}} )\,
\partial_{\underline{0}}.$ Since $T$ is local (cf. Prop. \ref{main6}) and
$\underline{\phi}$ is constant, it follows that (where $t_i$ for
$i=1,\ldots,l$ are smooth functions on $M$)
\begin{equation}\label{Eq1}
\left ( T( L_X ^{\underline{\l};\mu}A )\right
)(\underline{\phi})=\sum_{i=0}^l t_i\left( L_X
^{\mu-|\underline{\l}|}a_{\underline{0}}\right
)^{(i)}\underline{\phi}.
\end{equation}
On the other hand, \begin{equation} \label{Eq2} \left
(L_{X}^{\underline{\rho};\eta} T(A)\right
)(\underline{\phi})=L_{X}^{\eta-|\underline{\rho}|}\left (
\sum_{i=0}^l t_i\, a_{\underline{0}}^{(i)}\right
)\underline{\phi}.
\end{equation}
Since $T$ is an isomorphism, the function $t_0$ is not identically
zero. By comparing the coefficient of $t_0$ in Eq. \ref{Eq1} and
Eq.\ref{Eq2}, the result follows. \cqf\\
\section{The $\mathfrak{sl}(2)$-equivariant symbol
calculus}
Equivariant symbol calculus was carried out in \cite{bn} for the
case of binary differential operators. An explicit formula was
given for $k=2.$ In this section, we extend the results to $m$-ary
differential operators, exhibiting the symbol map for any $k.$ It
should be stressed, however, that the equivariant symbol calculus
depends on the embedding of the Lie algebra
$\mathfrak{sl}(2)\subset \Vect(M).$ For instance, if we consider
the embeddings
\[\mathfrak{sl}(2)\simeq \mbox{Span}\{\partial_x,
\sin(x)\partial_x,\cos(x)\partial_x\} \quad \mbox{or} \quad
\mathfrak{sl}(2)\simeq\mbox{Span}\{\partial_x,
\sinh(x)\partial_x,\cosh(x)\partial_x\},\] the symbol may not
exit, as already pointed out in \cite{by} for the unary case.
Throughout this paper, $\mathfrak{sl}(2)$ is realized as in
(\ref{Dub}).
\begin{thm} \label{main1} For all $\delta
\not\in\{1, \frac{3}{2},2,\ldots,k\},$ there exits a family of
$\mathfrak{sl}(2)$-equivariant maps given by
\[
\begin{array}{c}
\sigma_{\underline{\l},\mu}^\alpha:\displaystyle
\ccD_{\underline{\l};\mu}^k\longrightarrow \displaystyle
\bigoplus_{j=0}^{k} {\cal S}_{\delta-j}^{(j)}\quad A\mapsto
\displaystyle \sum_{r=0}^k\sum_{|\underline{i}|=r}
\overline{a}_{\underline{i}}(dx)^{\delta-|\underline{i}|}
\end{array}
\]
where $\overline{a}_{\underline{i}}=\displaystyle
\sum_{|\underline{s}|=|\underline{i}|}^k\alpha_{\underline{s}}^{\underline{i}}\,\,
a_{\underline{i}}^{(|\underline{s}|-|\underline{i}|)}$ and the
constants $\alpha^{\underline{s}}_{\underline{i}}$ are given by
the induction formula (where
$|\underline{s}|,|\underline{i}|=0,\ldots,k$):
\begin{equation}
\label{don} \left (|\underline{s}|-|\underline{i}|\right )\left
(2\delta-|\underline{s}|-|\underline{i}|-1\right
)\,\alpha^{\underline{s}}_{\underline{i}}-\sum_{j=1}^{m}
s_j(2\l_j+s_j-1)\,
\alpha_{\underline{i}}^{\underline{s}-\mathbf{1}_j}=0.
\end{equation}
\end{thm}
\noindent {\bf Proof.} Proposition \ref{mai} asserts that the
$\mathfrak{sl}(2)$-equivariant symbol map is local, hence it is
given by a differential operator. Let us study the invariance
property. We shall show that
\[
\overline{a}_{\underline{i}}^X=L^{\delta-|\underline{i}|}_{X}\,
\overline{a}_{\underline{i}}+ \mbox{ higher termes } X^{(n)},
n\geq 3.
\]
Now if we restrict ourself to $\mathfrak{sl}(2)$ then the second
part of the right hand side vanishes and thus we have
equivariance. To prove the formula above, we consider
\[
\overline{a}_{\underline{i}}=
\sum_{|\underline{s}|=|\underline{i}|}^k
\alpha_{\underline{i}}^{\underline{s}}\,
a_{\underline{s}}^{([\underline{s}]-[\underline{i}])}.
\]
Upon using proposition \ref{prop21} we get
\[
\small
\begin{array}{ccl}
\overline{a}_{\underline{i}}^X&=&L^{\delta-|\underline{i}|}_{X}\,
\overline{a}_{\underline{i}}\\[3mm]
& & \displaystyle +X''\sum_{|\underline{s}|=|\underline{i}|}^k
\alpha_{\underline{i}}^{\underline{s}} \left
(\binom{|\underline{s}-\underline{i}|}{2}+(\delta-|\underline{s}|)
|\underline{s}-\underline{i}|\right )
a_{\underline{s}}^{(|\underline{s}-\underline{i}|)}\\[3mm]
& &\displaystyle  -X''\,
\sum_{|\underline{s}|=|\underline{i}|-1}^{k-1}
\alpha_{\underline{i}}^{\underline{s}}\sum_{j=1}^m\,\left (\l_j
\binom{s_j+1}{1}+\binom{s_j+1}{2}\right )\,
a_{\underline{s}}^{(|\underline{s}-\underline{i}|+1)}+ \mbox{
higher terms in } X^{(n)}.
\end{array}
\]
The coefficient of $X''$ turns out to be trivial thanks to the
formula (\ref{don}) and the induction hypothesis. \cqfd

We have a family of symbol maps
$\sigma_{\underline{\l},\mu}^\alpha$ generated by the entries of
the non-singular matrices $[\alpha]_{i},$ for $i=1,\ldots,k.$
Nevertheless, we will prove that once the principal symbol is
fixed the symbol map $\sigma_{\underline{\l},\mu}^\alpha$ is
unique.
\begin{pro}
For $\delta\not \in \{1, \frac{3}{2}, 2,\ldots, k\},$ there exists
a unique $\mathfrak{sl}(2)$-equivariant symbol map
$\sigma_{\underline{\lambda};\mu}:{\cal
D}_{\underline{\lambda};\mu}\rightarrow {\cal S}_\delta$ such
that, for each $A\in {\cal D}_{\underline{\lambda};\mu}^k,$ the
term of highest order of $\sigma_{\lambda;\mu}$ is the symbole map
$\sigma^\alpha$ for some matrix $\alpha.$
\end{pro}
{\bf Proof.} The proof is similar to \cite{jgp}. Let assume that
there is another $\mathfrak{sl}(2)$-equivariant symbol map $\tilde
\sigma.$ Then, for every integer $k,$ the restriction of the map
$\sigma_{\underline{\lambda};\mu} \circ {\tilde \sigma}^{-1}$ to
${\cal S}_{\delta-k}^{(k)}$ is of the form
\[
a\in {\cal S}_{\delta-k}^{(k)} \mapsto
(S_0(a),S_1(a),\ldots,S_{k}(a))\in \bigoplus_{i=0}^k {\cal
S}_{\delta -i} ^{(i)}
\]
for some $\mathfrak{sl}(2)$-equivariant maps $S_i\in
\Hom_{\mathfrak{sl}(2)}({\cal S}_{\delta-k}^{(k)} ,{\cal
S}_{\delta-i}^{(i)} ),$ where $i=1,\ldots,k.$ We have
\[
\Hom_{\mathfrak{sl}(2)}({\cal S}_{\delta-k}^{(k)} ,{\cal
S}_\delta^{(i)} )= \bigoplus_{\binom{k+m-1}{m-1}
\mathrm{times}}\bigoplus_{\binom{i+m-1}{m-1}\mathrm{times}}
\Hom_{\mathfrak{sl}(2)}({\cal F}_{\delta-k} ,{\cal F}_{\delta-i}
).\] Following \cite{jgp} and since $\delta\not
\in\{1,\frac{3}{2}, 2,\ldots\}$ we have
\[
\Hom_{\mathfrak{sl}(2)}({\cal F}_{\delta-k} ,{\cal F}_{\delta-i}
)\simeq \left \{ \begin{matrix} \mathrm{Id}& \mathrm{if }\, i=k, \\
0& \mathrm{if }\, i\not =k.
\end{matrix}
\right.
\] Therefore, all the maps $S_i$ are zero except $S_0$ which is
given as a multiplication by a non-singular matrix since the maps
$\sigma_{\underline{\lambda};\mu}$ and $\tilde \sigma$ are
isomorphisms. The result follows.\cqf

We are interested in a class of symbol maps where the principal
symbols are given by
\[
[\alpha]_{i}=\mathrm{Id} \quad \mbox{for }\,i=1,\ldots,k.
\]
The inverse of the symbol map is the {\it quantization} map. It is
described by the following Theorem.
\begin{thm}
\label{main2} For all $\delta \not \in\{1,2,
\frac{3}{2},\ldots,k\}$ there exits a family of
$\mathfrak{sl}(2)$-equivariant maps given by
\begin{equation}
\label{main3}
\begin{array}{c}
\displaystyle Q_{\underline{\l},\mu}^\beta:\bigoplus_{j=0}^k{\cal
S}_{\delta-j}^{(j)} \longrightarrow
\ccD_{\underline{\l};\mu}^k\qquad \displaystyle
\bigoplus_{r=0}^k\bigoplus_{|\underline{i}|=r}
a_{\underline{i}}(dx)^{\delta-|\underline{i}|}\mapsto
\sum_{r=0}^{k}\sum_{|\underline{i}|=r}
\widetilde{a}_{\underline{i}}\,\partial_{\underline{i}}
\end{array}
\end{equation}
where $\widetilde{a}_{\underline{i}}=
\sum_{|\underline{s}|=|\underline{i}|}^k
\beta_{\underline{i}}^{\underline{s}}\,\,
a_{\underline{s}}^{(|\underline{s}-\underline{i}|)}$ and the
constants $\beta^{\underline{s}}_{\underline{i}}$ are given by the
induction formula (for
$|\underline{s}|,|\underline{i}|=0,\ldots,k):$
\begin{equation}
\label{Equi} \left (2\delta-1-|\underline{s}-\underline{i}|\right
)|\underline{s}-\underline{i}|\,\beta_{\underline{i}}^{\underline{s}}+\sum_{j=1}^m
(i_j+1)(2\l_j+i_j)\,\beta_{\underline{i}+\mathbf{1}_j}^{\underline{s}}=0.
\end{equation}
\end{thm}
{\bf Proof.} The proof is similar to that of Theorem \ref{main1}.
\cqfd
\section{Examples}
We provide examples of the quantization map (\ref{main3}) for the
case of first-order and second-order operators. These expressions
will be used to study the $\Vect(M)$-isomorphism problem.
\subsection{The case of first-order $m$-ary differential operators}
For first-order operators, the quantization map is given as in
(\ref{main3}), where the constants (\ref{Equi}) are given by (for
$j=1,\ldots,m$)
\begin{equation}
\label{Ex1} [\beta]_{1}=\mathrm{Id},\quad
\beta_{\underline{0}}^{\mathbf{1}_j}=
\displaystyle\frac{\l_j}{1-\delta}.
\end{equation}
\subsection{The case of second-order $m$-ary differential operators}
For second-order operators, the quantization map is given as in
(\ref{main3}), where the constants (\ref{Equi}) are given by (for
$i,j=1,\ldots,m$ and $i\not=j$)
$$
[\beta]_{2} =\mathrm{Id},\quad
\beta_{\underline{0}}^{2\mathbf{1}_j}=
\displaystyle\frac{\l_j(2\l_j+1)}{(\delta-2)(2\delta-3)},\quad
\beta_{\underline{0}}^{\mathbf{1}_i+\mathbf{1}_j}
=\displaystyle\frac{2\l_i\l_j} {(\delta-2)(2\delta-3)}.
$$
And
\begin{equation}
\label{Ex2}
\begin{array}{rclrclcrcl}
\beta_{\mathbf{1}_1}^{\mathbf{1}_1+\mathbf{1}_1}&=
&\displaystyle\frac{2\l_1+1}{2-\delta},&
\beta_{\mathbf{1}_1}^{\mathbf{1}_1+\mathbf{1}_2}&=
&\displaystyle\frac{\l_2}{2-\delta},&\cdots
&\beta_{\mathbf{1}_1}^{\mathbf{1}_1+\mathbf{1}_m}&=
&\displaystyle\frac{\l_m}{2-\delta},\\[3mm]
&\vdots &&&\vdots&&&&\vdots&\\[3mm]
\beta_{\mathbf{1}_m}^{\mathbf{1}_m+\mathbf{1}_m}&=
&\displaystyle\frac{2\l_m+1}{2-\delta},&
\beta_{\mathbf{1}_m}^{\mathbf{1}_m+\mathbf{1}_1}&=
&\displaystyle\frac{\l_1}{2-\delta},&\cdots
&\beta_{\mathbf{1}_m}^{\mathbf{1}_m+\mathbf{1}_{m-1}}
&=&\displaystyle\frac{\l_{m-1}}{2-\delta}.
\end{array}
\end{equation}
The other values of $\beta^{\underline{i}}_{\underline{s}}$
vanish.
\begin{rmk}{\rm
To another approach to the study of the space of $m$-ary
differential operators, see \cite{brov}. }
\end{rmk}
\section{The case of $m$-ary skew symmetric operators}
Consider now the space of $m$-ary skew symmetric differential
operators, $\ccD_{\wedge^m \l;\mu}.$ In that case, we deal only
with two parameters $\mu$ and $\l:=\l_1=\cdots=\l_m.$ Let
$Q(i,m)$\footnote{For the low values of $m,$ the function $Q(k,m)$
has beautiful expressions. For instance: (1) $Q(i,2)=\lfloor
\frac{1}{2}(i-1)\rfloor,$ where $\lfloor x \rfloor$ is the floor
function defined to be the greatest integer $\leq x$; (2)
$Q(i,3)=[ \frac{1}{12}(i-3)^2],$ where $[x]$ is the nint function
defined to be the closed integer to $x$ but half integers rounded
to even numbers, as $[1.5]=2,\, [2.5]=2.$ } be the number of ways
of partitioning $i$ into exactly $m$ distinct positive parts (see
\cite{com}). Let us put $R(i,m)=Q(i,m)+Q(i,m-1)$. The quotient
module $\ccD_{\wedge^m \l;\mu}^k/\ccD_{\wedge^m \l;\mu}^{k-1}$ can
be decomposed into $R(k,m)$-components that transform under
coordinate changes as $(\delta-k)$-densities. Denote by ${\cal
S}_{\l}^{R(i,m)}$ the direct sum $\oplus \cF_{\l}$ counted
$R(i,m)$-times. The symbol map is as follows:
\begin{thm} \label{main7} For all $\delta
\not\in\{1, \frac{3}{2},2,\ldots,k\},$ there exits a family of
$\mathfrak{sl}(2)$-equivariant maps given by
\[
\begin{array}{c}
\sigma_{\underline{\l},\mu}^\alpha:\displaystyle
\ccD_{\wedge^m\l;\mu}^k\longrightarrow \displaystyle
\bigoplus_{j=1}^{k} {\cal S}_{\delta-j}^{R(j,m)}\quad A\mapsto
\displaystyle
\sum_{r=1}^k\sum_{\substack{[\underline{i}]=r\\i_1>\cdots>
i_m}}\overline{a}_{\underline{i}}(dx)^{\delta-[\underline{i}]}
\end{array}
\]
where $\overline{a}_{\underline{i}}=\displaystyle
\sum_{|\underline{s}|=
|\underline{i}|}^{k}\alpha_{\underline{s}}^{\underline{i}}\,\,
a_{\underline{i}}^{(|\underline{s}-\underline{i}|)}$ and the
constants $\alpha^{\underline{s}}_{\underline{i}}$ are given by
the induction formula (for
$|\underline{s}|,|\underline{i}|=1,\ldots,k,$ where
$i_1>\cdots>i_m $ and $s_1>\cdots >i_m$):
\begin{equation}
\label{don9} |\underline{s}-\underline{i}| \left
(2\delta-|\underline{s}-\underline{i}|-1\right
)\,\alpha^{\underline{s}}_{\underline{i}}-\sum_{j=1}^{m}
s_j\,(2\l+s_j-1)\,
\alpha_{\underline{i}}^{\underline{s}-\mathbf{1}_j}=0.
\end{equation}
\end{thm}
{\bf Proof.} The operators that we are dealing with are skew
symmetric. Therefore, the components $a_{\underline{0}}$ must be
zero and the other components are skew symmetric with respect to
the indices. This explains why the sum should be taken over
distinct indices. Now the proof is similar to that of
Theorem \ref{main1}.\\
\cqfd
\begin{cor}
The following $\mathfrak{ sl}(2)$-modules are isomorphic
\[
\ccD^k_{\wedge^2\l;\mu+\l}\simeq \bigoplus_{i=0}^{\lfloor
\frac{1}{2}(k-1)\rfloor} \ccD^k_{\l;\mu}/\ccD^{i}_{\l;\mu}.
\]
\end{cor}
{\bf Proof.} This follows from Theorem \ref{main7} and the
equivariant quantization map exhibited in \cite{gar} for the unary
case. \\
\cqf
\begin{rmk}{\rm
Skew-symmetric invariant differential operators on weighted
densities have been investigated in \cite{ff}, generalizing the
Grozman operator \cite {groz} from $\ccD^3_{\wedge^2
\frac{-2}{3};\frac{5}{3}}.$ For a historical account, see
\cite{gls}.}
\end{rmk}
\section{The $\mathfrak{sl}(2)$-equivariant
quantization, the resonant case}
For the sake of completeness, we study the resonant values of
$\delta\in \{1, \frac{3}{2},2,\frac{5}{2},\ldots,k\}.$ The
following result contrasts with the unary case.
\begin{thm}
\label{main5} (i) For $\delta=1,$ a $\mathfrak{sl}(2)$-equivariant
map ${\cal S}_{\delta}\rightarrow \ccD_{\underline{\lambda};\mu}$
exists only for $\underline{\l}=\underline{0}.$

(ii) For $\delta=\frac{3}{2},$ a $\mathfrak{sl}(2)$-equivariant
map ${\cal S}_{\delta}\rightarrow \ccD_{\underline{\lambda};\mu}$
exists only for $\underline{\l}=\underline{0}$ or
$\underline{\l}=-\frac{1}{2}\mathbf{1}_j,$ where $j=1,\ldots,m.$

(iii) For $\delta\in \{2,\frac{5}{2},3,\ldots,k\}$ there is no
$\mathfrak{sl}(2)$-equivariant map ${\cal S}_{\delta}\rightarrow
\ccD_{\underline{\lambda};\mu}$, for any $\underline{\l}$ and
$\mu.$
\end{thm}
{\bf Proof.} The proof is based on a mathematical induction. The
$\mathfrak{sl}(2)$-equivariance is equivalent to the following
linear system (for $|\underline{s}-\underline{i}|\geq 0$ and
$|\underline{s}|,|\underline{i}|=0,\ldots,k$):
\begin{equation}
\label{out}
(2\delta-1-|\underline{s}-\underline{i}|)|\underline{s}-\underline{i}|
\,\beta_{\underline{i}}^{\underline{s}}+\sum_{j=1}^m
(i_j+1)(2\l_j+i_j)\,\beta_{\underline{i}+\mathbf{1}_j}^{\underline{s}}=0.
\end{equation}
In the case where $|\underline{s}|\leq k-1,$ the system
(\ref{out}) is exactly the $\mathfrak{sl}(2)$-equivariance
condition for the module $\ccD_{\underline{\l};\mu}^{k-1}.$ As we
have a filtration of modules
\[
\ccD_{\underline{\l};\mu}^1\subset
\ccD_{\underline{\l};\mu}^2\subset \cdots \subset
\ccD_{\underline{\l};\mu}^{k-1}\subset
\ccD_{\underline{\l};\mu}^{k},
\]
we will be dealing with the induction assumption at $k-1$ together
with the system (\ref{out}) for $|\underline{s}|=k.$

(i) The case where $\delta=1.$ Let us first study the case where
$k=1.$ The $\mathfrak{sl}(2)$-equivariance is equivalent to the
system
\[
\sum_{j=1}^m\beta^{\mathbf{1}_1}_{\mathbf{1}_j}\l_{j}
=0,\ldots,\sum_{j=1}^m\beta^{\mathbf{1}_m}_{\mathbf{1}_j}\l_{j}
=0.
\]
As the matrix $[\beta]_{1}$ is non-singular, it follows that
$\underline{\l}$ must be $\underline{0}.$ Suppose that the result
holds true at $k-1.$ Now, at $k$ we are required to solve the
system (\ref{out}) only for $|\underline{s}|=k$ -- actually for
$|\underline{s}|<k$ solutions are guaranteed by the induction
assumption. For this value, the system (\ref{out}) becomes
\begin{equation}
\label{out2} (1-k-|\underline{i}|)(k-|\underline{i}|)
\,\beta_{\underline{i}}^{\underline{s}}+\sum_{j=1}^m
(i_j+1)(2\l_j+i_j)\,\beta_{\underline{i}+\mathbf{1}_j}^{\underline{s}}=0.
\end{equation}
As $1-k-|\underline{i}|$ is never zero for every
$|\underline{i}|\leq k-1,$ the constant
$\beta_{\underline{i}}^{\underline{s}}$ can be expressed in terms
of $\beta_{\underline{i}+\mathbf{1}_j}^{\underline{s}}.$ This
means that every constant $\beta_{\underline{i}}^{\underline{s}}$
can be expressed in terms of the matrix $[\beta]_{k}.$ We have no
conditions on that matrix except that it should be non-singular.

(ii) The case where $\delta=\frac{3}{2}.$ This value is not a
resonant value for the module $\ccD^1_{\underline{\l};\mu},$ hence
the $\mathfrak{sl}(2)$-equivariant map exists. Here we cannot
proceed directly by induction because this value is resonant for
the module $\ccD^2_{\underline{\l};\mu}.$ So we need to prove the
result for $k=2,$ then we proceed by induction. The
$\mathfrak{sl}(2)$-equivariant is equivalent to the system:
\[
\begin{array}{rlcl}
\displaystyle -\beta_{\underline{i}}^{\underline{s}}+\sum_{j=1}^m
(i_j+1)(2\l_j+i_j)\,
\beta_{\underline{i}+\mathbf{1}_j}^{\underline{s}}&=&0 & \mbox{
for } |\underline{i}|=1,\\
\displaystyle \sum_{j=1}^m (i_j+1)(2\l_j+i_j)\,
\beta_{\underline{i}+\mathbf{1}_j}^{\underline{s}}&=&0 & \mbox{
for } |\underline{i}|=0.
\end{array}
\]
By solving this system, we get (for
$|\underline{s}|,|\underline{i}|=2$):
\[
\sum_{u=1}^m\sum_{v=1}^m(i_u-\delta^{v}_u)(2\l_u+i_u-\delta^v_u-1)
i_v(2\l_v+i_v-1)\beta^{\underline{s}}_{\underline{i}}=0.
\]
As the matrix $[\beta]_2$ is non-singular, the left part
$\sum_{u=1}^m\sum_{v=1}^m(i_u-\delta^{v}_u)(2\l_u+i_u-\delta^v_u-1)
i_v(2\l_v+i_v-1)$ must be zero. By taking
$\underline{i}=2\mathbf{1}_r$ (for $r=1,\ldots,m$) we obtain the
system \begin{equation} \label{eqq1}
4\l_r(2\l_r+1)=0.\end{equation} By taking
$\underline{i}=\mathbf{1}_p+\mathbf{1}_q$ (for $p,q=1,\ldots,m),$
we obtain the system
\begin{equation}\label{eqq2}\l_p\,\l_q=0.\end{equation}
The system (\ref{eqq1}, \ref{eqq2}) admits roots, given as stated
in the Theorem. Suppose that the result holds true at $k-1.$ As
explained in Part (i), we will deal only with the system
(\ref{out}) for $|\underline{s}|=k.$ For this value, the system
(\ref{out}) becomes
\begin{equation}
\label{out1} (2-k-|\underline{i}|)(k-|\underline{i}|)
\,\beta_{\underline{i}}^{\underline{s}}+\sum_{j=1}^m
(i_j+1)(2\l_j+i_j)\,\beta_{\underline{i}+\mathbf{1}_j}^{\underline{s}}=0.
\end{equation}
As the quantity $(2-k-|\underline{i}|)$ is never zero, for every
$|\underline{i}|\leq k-1$ and $k>2,$ then the constant
$\beta_{\underline{i}}^{\underline{s}}$ can be expressed in terms
of $\beta_{\underline{i}+\mathbf{1}_j}^{\underline{s}}.$ This
means that every constant $\beta_{\underline{i}}^{\underline{s}}$
can be expressed in terms of the matrix $[\beta]_{k}.$ We have no
conditions on that matrix except being non-singular. The result
follows upon taking the restriction to the submodule
$\ccD^{k-1}_{\underline{\l},\mu}.$

(iii) The case where $\delta\in \{2,\frac{5}{2},\ldots,k\}.$ Let
us start by studying the case where $\delta$ is an integer. The
proof can be obtained for $k=2$ but we here omit the details.
Suppose that the result is true at $k-1.$ The inclusion
$\ccD_{\underline{\l},\mu}^{k-1}\subset
\ccD_{\underline{\l};\mu}^k $ implies that the
$\mathfrak{sl}(2)$-equivariant map does not exist for
$\delta=2,\ldots,k-1.$ Let us prove the result for $\delta=k.$ For
this value the system (\ref{out}), for $|\underline{s}|=k,$
becomes (where $|\underline{i}|=k$):
\begin{equation}
\label{eqq3} \sum_{j=1}^m
i_j\,(2\l_j+i_j-1)\,\beta_{\underline{i}}^{\underline{s}}=0.
\end{equation}
For $i=2\mathbf{1}_1+(k-2)\mathbf{1}_2,$ the system (\ref{eqq3})
can be rewritten as
\[
[\beta]_k\times \left [
\begin{array}{c}
0\\
\vdots
\\
2(2\l_1+1)\\
(k-2)(2\l_2+k-3)\\
0\\
\vdots \\
0
\end{array}
 \right ]=
  \left [
\begin{array}{c}
0\\
\vdots
\\
0\\
0\\
0\\
\vdots \\
0
\end{array}
 \right ].
\]
For $i=\mathbf{1}_1+(k-1)\mathbf{1}_2,$ the system (\ref{eqq3})
can be rewritten as
\[
[\beta]_k\times \left [
\begin{array}{c}
0\\
\vdots
\\
2\l_1\\
(k-1)(2\l_2+k-2)\\
0\\
\vdots \\
0
\end{array}
 \right ]=
  \left [
\begin{array}{c}
0\\
\vdots
\\
0\\
0\\
0\\
\vdots \\
0
\end{array}
 \right ].
\]
As the matrix $[\beta]_k$ is not singular, it follows that
$\l_1=-\frac{1}{2}$ and $\l_1=0$ which is absurd.

Let us study the case where $\delta=\frac{2l-1}{2}$ for
$l=3,\ldots,k.$ Here the computation can be checked for $k=3$ but
we omit the details. Suppose that the result is true at $k-1.$ The
inclusion $\ccD_{\underline{\l},\mu}^{k-1}\subset
\ccD_{\underline{\l};\mu}^k $ implies that the
$\mathfrak{sl}(2)$-equivariant map does not exist for
$\delta=\frac{2l-1}{2},$ where $l=3,\ldots,k-1.$ Let us prove the
result for $\delta=\frac{2k-1}{2}.$ This value is, actually, not a
resonant value for the module $\ccD^{k-1}_{\underline{\l},\mu}.$
Therefore, the system (\ref{out}) admits a solution for every
$|\underline{s}|<k.$ Let us study this system for
$|\underline{s}|=k.$ By solving the system (\ref{out}) for
$|\underline{i}|=k-1$ we get
\begin{equation}
\label{eqq5}
\beta_{\underline{i}}^{\underline{s}}=
\sum_{j=1}^m(i_j+1)(2\l_j+i_j)\,
\beta_{\underline{i}+\mathbf{1}_j}^{\underline{s}}.
\end{equation}
Now for $|\underline{i}|=k-2,$ the system (\ref{out}) becomes
\begin{equation}
\label{eqq6}
 \sum_{j=1}^m
i_j\,(2\l_j+i_j-1)\beta_{\underline{i}}^{\underline{s}}=0.
\end{equation}
Upon substituting Eq. (\ref{eqq5}) into Eq. (\ref{eqq6}) we get
(for $|\underline{s}|,|\underline{i}|=k$)
\[
\sum_{u=1}^m\sum_{v=1}^m(i_u-\delta^{v}_u)(2\l_u+i_u-\delta^v_u-1)
i_v(2\l_v+i_v-1)\beta^{\underline{s}}_{\underline{i}}=0.
\]
As the matrix $[\beta]_k$ is non-singular, the left part
$\sum_{u=1}^m\sum_{v=1}^m(i_u-\delta^{v}_u)(2\l_u+i_u-\delta^v_u-1)
i_v(2\l_v+i_v-1)$ must be zero. By taking
$\underline{i}=k\mathbf{1}_r$ (for $r=1,\ldots,m$) we obtain the
system \begin{equation} \label{eqq7}
(1-k)(2\l_r+k-2)(2\l_r+k-1)=0.\end{equation} By taking
$\underline{i}=(k-1)\mathbf{1}_p+\mathbf{1}_q$ (for
$p,q=1,\ldots,m,\mbox{ and } p\not= q)$ we obtain the system
\begin{equation}\label{eqq8}(k-1)(2\l_p+k-2)\,
((k-2)(2\l_p+k-3)+4\l_q)=0.\end{equation} By taking
$\underline{i}=(k-2)\mathbf{1}_p+\mathbf{1}_q+\mathbf{1}_r$ (for
$p,q,r=1,\ldots,m,\mbox{ and } p,q,r \mbox{ are distinct})$ we
obtain the system
\begin{equation}\label{eqq9}(k-2)[(2\l_p+k-3)((k-3)(2\l_p+k-4)+
4\,\l_q +4\,\l_r)+8\l_q\l_r] =0.\end{equation} We distinguish two
cases:

1) If $\l_i=\frac{2-k}{2}$ for all $i=1,\ldots,m.$ By substituting
$\lambda_i$ in Eq. (\ref{eqq9}) we get
\[
\frac{1}{k}(2k-3)(k-2).
\]
This outcome is never zero for all $k\geq 3.$

2) If there exists $i_0$ such that $(2\l_{i_0}+k-2)\not=0.$ Eq.
(\ref{eqq7}) implies that $\l_{i_0}=\frac{1-k}{2}.$ Now Eq.
(\ref{eqq8}) implies that $\l_j=\frac{k-2}{2}$ for all $j\not=
i_0.$ By substituting in Eq. (\ref{eqq9}) we get \[
-2(k-2)(k-1)(2k-3).
\]
This last outcomes is never zero for all $k\geq 3.$

Thus, the system (\ref{out}) has no solutions and a fortiori there
is no $\mathfrak{sl}(2)$-equivariant quantization
map.\\
\cqfd
\section{A remark on $\Vect(M)$-equivariant quantization}
The $\mathfrak{sl}(2)$-equivariant quantization map is not unique,
generated by the entries of the matrices $[\beta]_i,$ where $
i=1,\ldots,k$. We can ask whether there exists an appropriate
principal symbol for which the equivariant quantization maps turn
into $\Vect(M)$-equivariant ones.
\begin{thm}
\label{main4} For $\delta\not \in \{1,\frac{3}{2},\ldots,k\},$
there exists a principal symbol for which the corresponding
quantization map is $\Vect(M)$-equivariant only in the following
cases:
\begin{enumerate}
\item For $k=1.$ \item For $k=2$ but
$\underline{\l}=\underline{0}$ or
$\underline{\l}=(1-\delta)\mathbf{1}_j$ for $j=1,\ldots,m.$
\end{enumerate}
\end{thm}
{\bf Proof.} We will first prove the result for $k=1, 2$ and $3.$
For $k=3,$ we will prove that no such principal symbol exists. As
we have a filtration of modules
\[{\cal D}^2_{\underline{\l},\mu}\subset
{\cal D}^3_{\underline{\l},\mu}\subset \cdots \subset {\cal
D}^k_{\underline{\l};\mu},\] the result holds for any $k>3$ upon
taking the restriction to the module ${\cal
D}^3_{\underline{\l},\mu}$ and applying the result.

For $k=1,$ the $\Vect(M)$-equivariance is given by the system
(\ref{Equi}) (for $|\underline{s}|,|\underline{i}|=0,1$). Upon
solving this system we get
\[
\beta_{\underline{0}}^{\mathbf{1}_j}=\sum_{s=1}^m\frac{\lambda_s}{1-\delta}
\beta_{\mathbf{1}_s}^{\mathbf{1}_j}\quad \mbox{for } j=1,\ldots,m.
\]
There are no more conditions on the constants
$\beta_{\mathbf{1}_s}^{\mathbf{1}_j}$ except that
$\mathrm{Det}[\beta]_1\not=0.$ We can take, for instance,
$[\beta]_1=\mathrm{Id},$ and therefore the corresponding
quantization map is certainly $\Vect(M)$-equivariant.

For $k=2,$ the $\Vect(M)$-equivariance is given by the system
(\ref{Equi}) (for $\underline{s},\underline{i}=0,1,2$) together
with the following system (for $\underline{u}=k$):
\begin{equation}
\label{eqq10}
 \sum_{j=1}^m \l_j
\beta^{\underline{u}}_{\mathbf{1}_j+\mathbf{1}_j}+(\delta-2)
\beta^{\underline{u}}_{\underline{0}}=0.
\end{equation}
By solving the system (\ref{Equi}), we get
\[
\beta_{\underline{0}}^{\underline{u}}=
\sum_{j=1}^m\frac{\l_j(2\l_j+1)}{(\delta-2)
(2\delta-3)}\beta^{\underline{u}}_{\mathbf{1}_j+\mathbf{1}_j}
+\sum_{\substack{i,j=1\\i\not= j}}^m\frac{2\l_j\l_i}{(\delta-2)
(2\delta-3)}\beta^{\underline{u}}_{\mathbf{1}_i+\mathbf{1}_j}.
\]
By substituting into Eq. (\ref{eqq10}) we get a new system (for
$|\underline{u}|=m$):
\[
\sum_{j=1}^m\lambda_j(\lambda_j+\delta-1)
\beta^{\underline{u}}_{2\mathbf{1}_j}
+\sum_{\substack{i,j=1\\i\not= j}}^m\l_i\l_j
\beta^{\underline{u}}_{\mathbf{1}_i+\mathbf{1}_j}=0.
\]
This system admits a solution for which the matrix $[\beta]_2$ is
non-singular if and only if the weights $\underline{\l}$ are given
as in Theorem \ref{main4}.

For $k=3,$ we proceed as above; two systems will be obtained that
we solve for the particular values of the weight $\underline{\l}.$
We omit details here but the proof is just a direct computation.\\
\cqfd
\section{Conjugation of $m$-ary differential operators}
First, we define a natural $\Vect(M)$-isomorphism on the modules
$\ccD_{\underline{\l};\mu}$ by just permuting arguments. Consider
the map $\mbox{per}_{i,j}$ that interchanges an element at the
$i$th position with an element at the $j$th position. This map
induces an isomorphism (for $i,j=1,\ldots,m \mbox{ and }
i~\not=~j~$)
\[
\ccD_{\underline{\l};\mu} \rightarrow
\ccD_{\mathrm{per}_{i,j}(\underline{\l});\mu} \quad A\mapsto
A\circ \mbox{per}_{i,j}.
\]
Now, we will define the notion of conjugation for $m$-ary
differential operators. For $M=\mathbb{R},$ We consider
compactly-supported densities.

Upon using successive integration by part, we get
\[
\int_{M} A(\varphi_1,\ldots,\varphi_m)\phi=\int_{M} \varphi_1
A^{*}(\varphi_2,\ldots,\varphi_m,\phi),
\]
where $A^{*}(\varphi_2,\ldots,\varphi_m,\phi)=\displaystyle
\sum_{\underline{i}} (-1)^{i_1}
\partial_{i_1}\,(a_{\underline{i}}\,\partial_{i_2}(\varphi_2)\,
\cdots \partial_{i_m}(\varphi_m)\,\phi).$ Therefore, the map $*$
induces a $\Vect(M)$-isomorphism
\[
\ccD_{\l_1,\ldots, \l_m;\mu}\xrightarrow{\simeq}
\ccD_{\l_2,\ldots,\l_m, 1-\mu;1-\l_1}\quad A\mapsto A^*.
\]
The following definition is adapted from the unary case
\cite{gar}.
\begin{defi}{\rm
A module $\ccD_{\underline{\l};\mu}$ is said to be singular if
either it is only isomorphic to itself, or it is isomorphic to any
another module ${\cal D}_{\underline{\rho};\eta}$ only through
compositions of conjugations and permutations. }
\end{defi}
\section{Classification of the modules $\ccD_{\underline{\l};\mu}^2$}
In this section we tackle the isomorphism problem. We study only
the case of second-order differential operators. The case where
$k>2$ seems to be more intricate. We need the following
\begin{pro} \label{aid} Every isomorphism
$T:\ccD_{\underline{\l};\mu}^k\rightarrow \ccD_{\underline{\rho};
\varrho}^k$ is block diagonal in terms of the
$\mathfrak{sl}(2)$-equivariant symbols. Namely, the map
$\sigma^{\mathrm{Id}}_{\underline{\l};\mu} \circ T\circ
Q^{\mathrm{Id}}_{\underline{\rho};\varrho}: \bigoplus_{i=0}^k{\cal
S}_{\delta-i}^{(i)}\rightarrow \bigoplus_{i=0}^k{\cal
S}_{\delta-i}^{(i)}$ is given by (where
$\tau^{\underline{i}}_{\underline{s}}$ are constants)
\begin{equation}
\label{eqq20}
\begin{array}{ccl}
\displaystyle
\bigoplus_{|\underline{i}|=0}^{k}\,a_{\underline{i}}&\mapsto&
\displaystyle \bigoplus_{|\underline{s}|=0}^k\,
\sum_{|\underline{i}|=|\underline{s}|}\,
\tau_{\underline{s}}^{\underline{i}}\,\, a_{\underline{i}},
\end{array}
\end{equation}
and $[\tau]_{i}$ are non-singular matrices for $i=0,1,\ldots,k.$
\end{pro}
{\bf Proof.} As $T$ is $\Vect(M)$-equivariant, it follows that the
composition
\begin{equation}
\ccD_{\underline{\l};\mu}^k
\xrightarrow{T}\ccD_{\underline{\rho};\eta}^k
\xrightarrow{\sigma_{\underline{\rho};\eta}^{\mathrm{Id}}}\bigoplus_{i=0}^k
{\cal S}_{\delta-i}^{(i)}
\end{equation}
is $\mathfrak{sl}(2)$-equivariant. Therefore, it coincides with
the symbol map $\sigma_{\underline{\l}; \mu}^{\tau}$ for some
$\tau.$ Namely, $\sigma^{\mathrm{Id}}_{\underline{\rho};\eta}\circ
T=\sigma_{\underline{\l}; \mu}^{\tau}.$ It follows that
\[
\sigma^{\mathrm{Id}}_{\underline{\rho};\eta}\circ T \circ
Q^{\mathrm{Id}}_{\underline{\l};\mu}=\sigma_{\underline{\l};
\mu}^{\tau} \circ Q^{\mathrm{Id}}_{\underline{\l};\mu}.
\]
Now, it is a matter of a direct computation to prove that
$\sigma_{\underline{\l}; \mu}^{\tau} \circ
Q^{\mathrm{Id}}_{\underline{\l};\mu}$ is given as (\ref{eqq20}).\\
\cqf
\subsection{The generic case}
We start by studying the case where $\delta\not = 1, \frac{3}{2},
2.$
\begin{thm}
(i) For $\delta\not=1,$ all modules $\ccD_{\underline{\l};\mu}^1$
are isomorphic provided they have the same shift $\delta.$\\
(ii) For $\delta\not =1,\frac{3}{2},2,$ all modules
$\ccD_{\underline{\l};\mu}^2$ are isomorphic provided they have
the same shift $\delta;$ however, the modules
\[
\ccD_{\underline{0};\delta}^2\simeq
\ccD_{(1-\delta)\mathbf{1}_1;1}^2\simeq\ldots
\simeq \ccD_{(1-\delta)\mathbf{1}_m;1}^2
\]
are singular.
\end{thm}
{\bf Proof.} By virtue of Proposition \ref{aid} we deal with the
map
\[
\sigma^{\mathrm{Id}}_{\underline{\l};\mu} \circ T\circ
Q^{\mathrm{Id}}_{\underline{\rho};\varrho}: \bigoplus_{i=0}^k{\cal
S}_{\delta-i}^{(i)}\rightarrow \bigoplus_{i=0}^k{\cal
S}_{\delta-i}^{(i)}\quad
\bigoplus_{|\underline{i}|=0}^{k}\,a_{\underline{i}}\mapsto
\displaystyle \bigoplus_{|\underline{s}|=0}^k\,
\sum_{|\underline{i}|=|\underline{s}|}\,
\tau_{\underline{s}}^{\underline{i}}\,\, a_{\underline{i}}.
\]
We are required to exhibit the coefficients
$(\tau_{\underline{i}}^{\underline{s}}).$ For this matter, we need
to compute the action
\begin{equation}
\label{eqq12} \sigma^{\mathrm{Id}}_{\underline{\l};\mu}\circ
L_{X}^{\underline{\l},\mu} \circ
Q^{\mathrm{Id}}_{\underline{\rho};\varrho}.
\end{equation}
For Part (i), a direct computation shows that the action on ${\cal
S}_{\delta-1}^{(1)}\oplus {\cal S}_\delta$ reads as follows
$$
\begin{array}{ccl}
\overline{a}^X_{\underline{i}}&=&L_{X}
^{\delta-1}a_{\underline{i}} \quad
\mbox{ for } |\underline{i}|=1,\\[2mm]
\overline{a}^X_{\underline{0}}&=&L_{X} ^{\delta}a_{\underline{0}}.
\end{array}
$$
Therefore, the modules $\ccD^1_{\underline{\l};\mu}$ and ${\cal
S}_{\delta-1}^{(1)}\oplus {\cal S}_\delta$ are isomorphic to each
other. We can choose the parameters
$\tau^{\underline{s}}_{\underline{i}}$ as
\[
[\tau]_{1}=\mathrm{Id} \quad \mbox{ and } \tau_{\underline{0}}=1.
\]
For Part (ii), a direct computation shows that the action
(\ref{eqq12}) on ${\cal S}_{\delta-2}^{(2)}\oplus {\cal
S}_{\delta-1}^{(1)}\oplus {\cal S}_\delta$ reads as follows:
\begin{equation}
\label{eqq21}
\begin{array}{ccl}
\overline{a}^X_{\underline{i}}&=&L_{X}
^{\delta-2}a_{\underline{i}} \quad
\mbox{ for } |\underline{i}|=2,\\[2mm]
\overline{a}^X_{\underline{i}} &=&L_{X}
^{\delta-1}a_{\underline{i}}
\quad \mbox{ for } |\underline{i}|=1,\\[2mm]
\overline{a}^X_{\underline{0}}&=&L_{X}
^{\delta}a_{\underline{0}}+\displaystyle \sum_{|\underline{i}|=2}
\alpha_{\underline{i}}\, X'''\, a_{\underline{i}},
\end{array}
\end{equation}
where (for $s,t=1,\ldots,m$ and $s\not =t$ ):
\begin{equation}\nonumber \alpha_{2\mathbf{1}_s}=2\l_s
\frac{1-\delta-\l_s}{2\delta-3}, \quad
\alpha_{\mathbf{1}_s+\mathbf{1}_t}=-\displaystyle 2 \frac{\l_s
\lambda_t}{2\delta-3}.
\end{equation}
The action (\ref{eqq21}) cannot be the action (\ref{act}) because
the 1-cocycle
\[
\Vect(M)\rightarrow \ccD_{\theta;\theta+2}\quad (X,\phi)\mapsto
X'''\phi,
\]
is not trivial for $\theta\not=-\frac{1}{2}$ (cf. \cite{bo,f}). We
define the column matrix (of $\binom{m+1}{m-1}$-entries) by
\[
\alpha(\underline{\l},\mu)=\left [
\begin{array}{c}
\alpha_{2\mathbf{1}_1}(\underline{\l},\mu)\\
\alpha_{\mathbf{1}_1+\mathbf{1}_2}(\underline{\l},\mu)\\
\vdots\\
\alpha_{2\mathbf{1}_m}(\underline{\l},\mu)
\end{array}
\right ].
\]
The existence of the $\Vect(M)$-isomorphism is equivalent to
solving the linear system
\begin{equation}
\label{linco} [\tau]_{2}\cdot
\alpha(\underline{\l},\mu)=\alpha^t(\underline{\rho}, \varrho)
\quad \mbox { and }\quad \mathrm{Det} [\tau]_{2}\not =0.
\end{equation}
We distinguish two cases:

1) If all entries of the column matrix
$\alpha(\underline{\l},\mu)$ are zero so are the entries of the
row matrix $c(\underline{\rho},\varrho),$ as $\mathrm{Det}
[\tau]_{2}\not =0.$ Now, the roots of the equation
$\alpha(\underline{\l},\mu)=0$ are
\begin{equation}
\label{sing} \underline{\l}=\underline{0}\quad \mbox{ or }\quad
(1-\delta)\mathbf{1}_j \quad \mbox{for } j=1,\ldots,m.
\end{equation}
Besides, the values of $\underline{\rho}$ must also be in the form
(\ref{sing}). Therefore, the corresponding isomorphism is the
composition of permutations and conjugations. Thus, the modules
(where $\underline{\l}$ is as in (\ref{sing}))
\[
\ccD_{\underline{\l};\mu}^2
\]
are singular.

2) If the column matrix $\alpha(\underline{\l},\mu)$ is not
identically zero neither is the column matrix
$\alpha(\underline{\rho},\varrho).$ Whatever the weights
$\underline{\l}$ and $\underline{\rho}$ are, the constant
$\tau^{\underline{i}}_{\underline{s}}$ can be chosen such that the
conditions (\ref{linco}) are satisfied. Thus, all modules
$\ccD_{\underline{\l};\mu}^2$ are isomorphic to each other. \cqfd
\begin{rmk}{\rm
The non-uniqueness of the isomorphism $T$ is also a characteristic
feature of the $m$-ary case. }
\end{rmk}
\subsection{The resonant case; the case of binary operators}
Throughout this section we deal with $m=2;$ thus for instance
$\underline{\l}$ stands for $(\l_1,\l_2).$ We shall study the case
when $\delta=1,\frac{3}{2},2.$ The quantization map exists only
for some particular values of $\underline{\l}.$ Hence the
techniques used in the previous section do not work. Here, we
proceed explicitly.
\begin{thm}
For $\delta=1,$ all modules
$\ccD_{\underline{\l};1+[\underline{\l}]}^1$ are isomorphic.
However, we have one exceptional module
\[
\ccD_{\underline{0};1}^1.
\]
\end{thm}
{\bf Proof.} We establish an isomorphism between the modules
$\ccD^1_{\underline{\l};[\underline{\l}]+1} \rightarrow
\ccD^1_{\underline{\rho};[\underline{\rho}]+1}$ as follows

1) If $\rho_2=0,$ then take $T$ as
\[
(a_{\underline{i}},a_{\underline{0}})\mapsto \left
(\sum_{[\underline{s}]=1}\tau_{\underline{i}}^{\underline{s}}\,
a_{\underline{s}},a_{\underline{0}}\right )\quad \mbox{for} \quad
[\underline{i}]=1,
\]
where the matrix $[\tau]_{1}=\left [
\begin{array}{cc}
\frac{\l_1}{\rho_1} & \frac{\l_2}{\rho_2}\\
c_1 & c_2
\end{array}
\right ].$ Whatever the values $\l_1$ and $\l_2$ can take, the
constants $c_1$ and $c_2$ can be chosen in a way such that
$\mathrm{Det}[\tau]_{1}\not=0.$

2) If $\rho_2,\rho_1\not=0,$ then take $T$ as
\[
(a_{\underline{i}},a_{\underline{0}})\mapsto \left
(\sum_{[\underline{s}]=1}\tau_{\underline{i}}^{\underline{s}}\,
a_{\underline{s}},a_{\underline{0}}\right )\quad \mbox{for} \quad
[\underline{i}]=1,
\]
where the matrix $[\tau]_{1}=\left [
\begin{array}{cc}
\frac{\l_1}{\rho_1} & 0\\
0 & \frac{\l_2}{\rho_2}
\end{array}
\right ].$ We point out that $\l_1, \l_2 \not =0;$ otherwise, we
go back to Part 1.

Suppose now that $\ccD^1_{\underline{\l};[\underline{\l}]+1}$ is
isomorphic to $\ccD^1_{0;1}.$ Therefore, the composition map
\[
\ccD^1_{\underline{\l};[\underline{\l}]+1}\rightarrow \ccD^1_{0;1}
\rightarrow \cF^{(1)}_{\delta-1}\oplus \cF_{\delta},
\]
is a $\mathfrak{sl}(2)$-equivariant quantization map in
contradiction with Theorem \ref{main5}. \cqfd

\begin{thm}
For $\delta=2,$ we have two classes of binary differential operators
\[
\ccD_{s\mathbf{1}_2;2+s}^2 \simeq \ccD_{s\mathbf{1}_1;2+s}^2 \quad
\mbox{ and }\quad \ccD_{\underline{\l};2+[\underline{\l}]}^2
\,\,(\l_1,\l_2 \not=0).
\]
However, we have two exceptional modules:
\[
\ccD_{-\mathbf{1}_1;1}^2 \simeq \ccD_{-\mathbf{1}_2;1}^2\simeq
\ccD_{\underline{0};2}^2\quad \mbox{ and } \quad
\ccD_{-\frac{1}{2}\mathbf{1}_1;\frac{3}{2}}^2
\simeq\ccD_{-\frac{1}{2}\mathbf{1}_2;\frac{3}{2}}^2\simeq
\ccD_{\underline{-1/2};1}^2\quad \mbox{(conjugations)}.\] (iii)
For $\delta=1,$ all modules are isomorphic. However, we have one
singular module
\[
\ccD_{\underline{0};1}^2.
\]
(iv) For $\delta=3/2,$ all modules are isomorphic. However we have
one singular module
\[
\ccD^2_{-\frac{1}{2}\mathbf{1}_1;1} \simeq
\ccD^2_{-\frac{1}{2}\mathbf{1}_2;1} \simeq
\ccD^2_{\underline{0};\frac{3}{2}} \mbox{ (conjugations) }.
\]
\end{thm}
{\bf Proof.} By using Proposition \ref{main6}, every isomorphism
$T:\ccD^2_{\underline{\l};\mu} \rightarrow
\ccD^2_{\underline{\rho};\eta}$ is local. Therefore, the map $T$
retains the following general form (where $[i]=2$ and $[j]=1$)
\[
\begin{array}{l}
(a_{\underline{i}},a_{\underline{j}},a_{\underline{0}})\mapsto\\[3mm]
\displaystyle \left
(\sum_{[\underline{s}]=2}\tau_{\underline{i}}^{\underline{s}}\,a_{\underline{s}},
\sum_{[\underline{s}]=2} \tau_{\underline{j}}^{\underline{s}}\,
a_{\underline{s}}'+\sum_{[\underline{s}]=1}
\tau^{\underline{s}}_{\underline{j}}\, a_{\underline{s}},
\sum_{[\underline{s}]=2} \tau^{\underline{s}}_{\underline{0}}\,
a_{\underline{s}}''+\sum_{[\underline{s}]=1}
\tau_{\underline{0}}^{\underline{s}}\,a_{\underline{s}}'
+\tau_{\underline{0}}\, a_{\underline{0}}\right ).
\end{array}
\]
A long and tedious computation proves that the
$\Vect(\mathbb{R})$-equivariant property is equivalent to the
following system of fourteen equations:
\begin{eqnarray*}
\lambda_1\,\tau^{\mathbf{1}_2}_{\underline{0}}
+\lambda_2\,\tau^{\mathbf{1}_1}_{\underline{0}}
-(\delta-2)\tau^{\mathbf{1}_1+\mathbf{1}_2}_{\underline{0}}-\sum_{j=1}^2\rho_j\,
\tau^{\mathbf{1}_1+\mathbf{1}_2}_{2\mathbf{1}_j}&=&0\\
(1+2\l_s)\,\tau^{\mathbf{1}_s}_{\underline{0}}-(2\delta-3)
\tau^{2\mathbf{1}_s}_{\underline{0}}-\sum_{j=1}^2\rho_j\,
\tau^{2\mathbf{1}_s}_{\mathbf{1}_j}&=&0\quad \mbox{ for } s=1,2\\
\l_1\,\tau^{\mathbf{1}_2}_{\underline{0}}+
\l_2\,\tau^{\mathbf{1}_1}_{\underline{0}}
-(2\delta-3)\tau^{\mathbf{1}_1+\mathbf{1}_2}_{\underline{0}}-
\sum_{j=1}^2\rho_j\,
\tau^{\mathbf{1}_1+\mathbf{1}_2}_{\mathbf{1}_i}&=&0\\
\l_s-(\delta-1)\,
\tau^{\mathbf{1}_s}_{\underline{0}}-\sum_{j=1}^2\rho_j\,
\tau^{\mathbf{1}_s}_{\mathbf{1}_j}&=&0\quad \mbox{ for } s=1,2\\
(1+2\l_s)\,\tau^{\mathbf{1}_s}_{\mathbf{1}_j}-(\delta-2)\tau^{2\mathbf{1}_s}_{\mathbf{1}_j}
-(1+2\rho_s)\,\tau^{2\mathbf{1}_s}_{2\mathbf{1}_j} -\rho_j\,
\tau^{2\mathbf{1}_s}_{\mathbf{1}_1+\mathbf{1}_2}&=&0 \quad \mbox{
for } s,j=1,2\\
\sum_{i=1}^2\l_j\,\tau^{\mathbf{1}_j}_{\mathbf{1}_s}
-(\delta-2)\,\tau^{\mathbf{1}_1+\mathbf{1}_2}_{\mathbf{1}_s}-
(1+2\rho_s)\, \tau^{\mathbf{1}_1+\mathbf{1}_2}_{2\mathbf{1}_s}
-\rho_{s+1}\,\tau^{\mathbf{1}_1+\mathbf{1}_2}_{\mathbf{1}_1+\mathbf{1}_2}&=&0
\quad \mbox{ for } s=1,2\\
\lambda_s+(1+2\l_s)\,\tau^{\mathbf{1}_s}_{\underline{0}}
-(\delta-2)\tau^{2\mathbf{1}_s}_{\underline{0}}-\sum_{i=1}^2\rho_i\,
\tau^{2\mathbf{1}_s}_{2\mathbf{1}_i}&=&0\quad \mbox{ for } s=1,2
\end{eqnarray*}
We give the details of the computation only for
$\delta=\frac{3}{2}.$ Here we distinguish also many cases:

1) If $\rho_1,\rho_2\not=0,$ then
\[
[\tau]_{1}=\mathrm{Id}\quad \mbox{ and } \quad [\tau]_{2}= \left[
\begin{array}{ccc}
  1 & \frac{(\l_1-\rho_1)(1+2\l_1+2\rho_1)}{2\rho_1 \rho_2} & 0
  \\[2mm]
  0 & \frac{(\l_2-\rho_2)(1+2\l_2+2\rho_2)}{2\rho_1 \rho_2} & 1
  \\[2mm]
  0 & \frac{\l_1\l_2}{\rho_1\rho_2} & 0 \\
\end{array}
\right],
\]
together with (where $s,j=1,\ldots,m$)
\[
\begin{array}{rclrcl}
\tau^{2\mathbf{1}_s}_{\underline{0}}&=&(6+8\l_s)(\rho_s-\l_s),&
\tau^{\mathbf{1}_1+\mathbf{1}_2}_{\underline{0}}&
=&4(\l_2\rho_1 +\lambda_1(\rho_2-2\l_2)),\\[2mm]
\tau^{\mathbf{1}_s}_{\underline{0}}&=& 2(\l_s-\rho_s),&
\tau^{\mathbf{1}_1+\mathbf{1}_2}_{\mathbf{1}_s}&=&
2\l_{s+1}(-1+\frac{\lambda_s}{\rho_s}),\\[2mm]
\tau^{2\mathbf{1}_s}_{\mathbf{1}_j}&=&\frac{(\lambda_s-\rho_s)
(1+2\lambda_s-2\rho_s)}{\rho_j}.
\end{array}
\]
2) If $\rho_1=0$ but $\rho_2\not=0,$ then
\[
[\tau]_{1}=\mathrm{Id}\quad \mbox{ and } \quad [\tau]_{2}= \left[
\begin{array}{ccc}
  1 & 0 & 0
  \\[2mm]
  0 & x & y
  \\[2mm]
  \frac{\l_1(1+2\l_1)}{\rho_2(1+2\rho_2)}&
  \frac{2\l_1\l_2}{\rho_2(1+2\rho_2)} &
  \frac{\l_2(1+2\l_2)}{\rho_2(1+2\rho_2)}  \\
\end{array}
\right],
\]
where the constants $x$ and $y$ can be chosen in a way such that
$\mbox{Det}[\tau]_2\not=0.$ The other constants are given by
\[
\begin{array}{rclrcl}
\tau^{2\mathbf{1}_2}_{\underline{0}}&=&-\frac{4(\l_2-\rho_2)
(1+\l_2+2\rho_2+4\l_2\rho_2)}{(1+2\rho_2)},&
\tau^{2\mathbf{1}_1}_{\underline{0}}&=&-\frac{4\l_1
(1+\l_1+3\rho_2+4\l_1\rho_2)}{(1+2\rho_2)},\\[2mm]
\tau^{\mathbf{1}_2}_{\underline{0}}&=&
2(\l_2-\rho_2),&\tau^{\mathbf{1}_1}_{\underline{0}}&=& 2\l_1,\\[2mm]
\tau^{2\mathbf{1}_2}_{\mathbf{1}_1}&=&2\rho_2 y, &
\tau^{\mathbf{1}_1+\mathbf{1}_2}_{\mathbf{1}_1}&=&
2(-\l_2+x\, \rho_2),\\[2mm]
\tau^{2\mathbf{1}_2}_{\mathbf{1}_1}&=&-4\l_1, &
\tau^{2\mathbf{1}_2}_{\mathbf{1}_2}&=&
\frac{2(1+2\l_2)(\l_2-\rho_2)}{\rho_2},\\[2mm]
\tau^{\mathbf{1}_1+\mathbf{1}_2}_{\mathbf{1}_2}&=&
\l_1(-2+\frac{4\l_2}{\rho_2}), &
\tau^{2\mathbf{1}_2}_{\mathbf{1}_1}&=&
\frac{2\l_1(1+2\l_1)}{\rho_2}.
\end{array}
\]

Theorem \ref{main4} asserts that
$\ccD^2_{\underline{\l};[\underline{\l}]+\frac{3}{2}},$ for
$\underline{\l}\not =
\underline{0},-\frac{1}{2}\mathbf{1}_1,-\frac{1}{2}\mathbf{1}_2,$
is not isomorphic to $\ccD^2_{\underline{0};\frac{3}{2}}\simeq
\ccD^2_{-\frac{1}{2}\mathbf{1}_1;1}\simeq
\ccD^2_{-\frac{1}{2}\mathbf{1}_2;1}.$ \cqfd
\subsection*{Acknowledgments}
Many thanks to V. Ovsienko for his encouragements and many thanks
to D. Leites, F. Wagemann and the reviewer for pointing out
pertinent remarks.


\begin{thebibliography}{99}
\pagestyle{myheadings} \markboth{\small Bibliographie}{\small
Bibliographie}

\bibitem{bor} M. Bordemann, Sur l'existence d'une prescription d'ordre naturelle
projectivement invariante, \verb"math.DG/0208171."

\bibitem{bl} S. Bouarroudj, Projectively equivariant quantization map,
{\it Lett. Math. Phys.,} {\bf 51,} no.4, (2000), 265--274.

\bibitem{by} S. Bouarroudj and M. Iyadh Ayari,
On $\mathfrak{sl}(2,\mathbb{R})$-equivariant quantizations, to
appear in {\it J. Nonlinear Math. Phys.} Available at
\verb"math.DG/0601353."

\bibitem{bo} S. Bouarroudj and V. Ovsienko, {\rm Three cocycles on
$\Diff(S^1)$ generalizing the Schwarzian derivative}, {\it
Internat. Math. Res. Notices,} no.1, (1998), 25--39.

\bibitem{bn} F. Boniver, {\rm  Projectively equivariant symbol
calculus for bidifferential operators,} {\it Lett. Math. Phys.}
{\bf 54,} no.2, (2000), 83-100.

\bibitem{bhm} F. Boniver, S. Hansoul, P. Mathonet and N. Poncin,
{\rm  Equivariant symbol calculus for differential operators
acting on forms,} {\it Lett. Math. Phys.} {\bf 62,} no.3, (2002),
219-232.

\bibitem{cmz} P. Cohen, Yu. Manin and D. Zagier, Automorphic
pseudo-differential operators,  in Algebraic Aspects of Integrable
Systems, Prog. Nonlinear Differential Equations Appl., {\bf 26},
Birkh\"{a}user, Boston, 1997, 17--47.

\bibitem{com} L. Comtet, {\it Advanced combinatorics: the art of finite
and infinite expansions,} rev. enl. ed. Dordrecht, Netherlands:
Reidel, 1974.

\bibitem{brov} V. Dobrev, {\rm
New generalized Verma modules and multilinear intertwining
differential operators,} {\it J.Geom.Phys.} {\bf 25} (1998) 1--28.

\bibitem{dlo}
C. Duval, P. Lecomte and V. Ovsienko, {\rm Conformally equivariant
quantization: existence and uniqueness.} {\it Ann. Inst. Fourier}
{\bf 49,} no. 6, (1999), 1999-2029.

\bibitem{do}
C. Duval and V. Ovsienko, {\rm Space of second order linear
differential operators as a module over the Lie algebra of vector
fields}, {\it Adv. in Math.} {\bf 132,} no. 2, (1997), 316--333.

\bibitem{ff} B. L. Feigin \& D. B. Fuchs, {\rm Invariant
skew-symmetric differential operators on the line and Verma
modules over the Virasoro algebra,} {\it Funkts. Anal.
Prilozhen.,} Vol. 16, no. 2, (1982), 47--63.

\bibitem{f} D. B. Fuks, {\it Cohomology of infinite-dimensional
Lie algebras}, Contemp. Soviet. Math., Consultants Bureau,
New-York, 1986.

\bibitem{gar}
H. Gargoubi, {\rm Sur la g\'eom\'etrie de l'espace des
op\'erateurs diff\'erentiels lin\'eaires sur $\mathbb{R}.$} {\it
Bull. Soc. Roy. Sci. Li\`ege.} Vol. 69, 1, (2000), 21--47.

\bibitem{groz}
P. Ya Grozman, {\rm Classification of bilinear invariant operators
over tensor fields,} {\it Functional Anal. Appl.,} {\bf 14}, no.2,
(1980), 127--128.

\bibitem{gls} P. Grozman, D. Leites, I. Shchepochkina, {\rm
Invariant operators on supermanifolds and standard models,} in:
Multiple facets of quantization and supersymmetry, 508-555, M.
Olshanetski, A. Vainstein (Eds.), Wolrd Sci. Publishing, 2002.

\bibitem{hans}
S. Hansoul, {\rm Projectively equivariant quantization for
differential operators acting on forms} {\it Lett. Math. Phys.}
{\bf 70,} no. 2, (2004), 141--153.

\bibitem{kir}
A. A. Kirillov, {\rm Invariant operators over geometric quantities
(Russian), in: Current Problems in Mathematics,} {\bf 16}, 3--29,
Akad. Nauk SSSR, VINITI, Moscow, 1980; [English translation: {\it
J. Sov. Math.} {\bf 18:1} (1982), 1--21].

\bibitem{jgp}
P. B. A. Lecomte, {\rm On the cohomology of
$\mathfrak{sl}(m+1,\mathbb{R})$ acting on differential operators
and $\mathfrak{sl}(m+1,\mathbb{R})$-equivariant symbol}, {\it
Indag. Math., N.S.,} {\bf 11,} no. 1, (2000), 95--114.

\bibitem{lmt}
P. B. A. Lecomte, P. Mathonet and E. Tousset, {\rm Comparison of
some modules of the Lie algebra of vector fields}, {\it Indag.
Math., N.S.,} {\bf 7,} no. 4, (1996), 461--471.

\bibitem{lo1}
P. B. A. Lecomte and V. Ovsienko, Projectively invariant symbol
calculus, {\it Lett. Math. Phy.,} {\bf 49,} no. 3, (1999),
173--196.

\bibitem{loub}
S. E. Loubon Djounga, Modules of third-order differential
operators on a conformally flat manifold. {\it J. Geom. Phys.,}
{\bf 37,} no. 3, (2001), 251--261.

\bibitem{mr}
P. Mathonet and F. Redoux, {\rm Projectively equivariant
quantizations by means of Cartan connections,} {\it Lett. Math.
Phys.,} {\bf 72,} no. 3, (2005), 183--196.

\bibitem{Pee} J. Petree, {\rm Une caract\'erisation abstraite des op\'erateurs
diff\'erentiels.} {\it Math. Scand.} {\bf 7} (1959), 211--218 and
{\bf 8}, (1960), 116--120.

\bibitem{wil} E. J. Wilczynski, {\it Projective differential geometry
of curves and ruled surfaces,} Leipzig, Teubner, 1906.


\end{thebibliography}
\end{document}